\documentclass[oneside,a4paper]{article}%
\usepackage[final]{microtype}%
\usepackage[colorlinks=true,linkcolor={blue!30!black}]{hyperref}%
\usepackage{amsthm,mathtools}%
\usepackage{autonum}%
\usepackage{thmtools}%
\usepackage{xcolor}%
\usepackage{float}%
\usepackage[export]{adjustbox}%
\usepackage{tikz, tikz-cd, mathtools, amssymb, stmaryrd}\usetikzlibrary{calc}
\usetikzlibrary{decorations.pathmorphing}
\usetikzlibrary{spath3}
\usetikzlibrary{backgrounds}
\usetikzlibrary{arrows}
\usetikzlibrary{shapes,shapes.geometric,shapes.misc}
\usepackage{circuitikz}

\tikzset{curve/.style={settings={#1},to path={(\tikztostart)
    .. controls ($(\tikztostart)!\pv{pos}!(\tikztotarget)!\pv{height}!270:(\tikztotarget)$)
    and ($(\tikztostart)!1-\pv{pos}!(\tikztotarget)!\pv{height}!270:(\tikztotarget)$)
    .. (\tikztotarget)\tikztonodes}},
    settings/.code={\tikzset{quiver/.cd,#1}
        \def\pv##1{\pgfkeysvalueof{/tikz/quiver/##1}}},
    quiver/.cd,pos/.initial=0.35,height/.initial=0}
\tikzset{between/.style n args={2}{/tikz/spath/at end path construction={
    \tikzset{spath/split at keep middle={current}{#1}{#2}}
}}}

\tikzset{tail reversed/.code={\pgfsetarrowsstart{tikzcd to}}}
\tikzset{2tail/.code={\pgfsetarrowsstart{Implies[reversed]}}}
\tikzset{2tail reversed/.code={\pgfsetarrowsstart{Implies}}}
\tikzset{no body/.style={/tikz/dash pattern=on 0 off 1mm}}
    \usepackage{mathpartir}
\usepackage{newpxmath,newpxtext}%
\usepackage{cleveref}%
\newtheorem{theorem}{Theorem}[section]%
\newtheorem{lemma}[theorem]{Lemma}%
\newtheorem{corollary}[theorem]{Corollary}%
\theoremstyle{definition}%
\newtheorem{definition}[theorem]{Definition}%
\newtheorem{proposition}[theorem]{Proposition}%
\newtheorem{example}[theorem]{Example}%
\newtheorem{remark}[theorem]{Remark}%
\newtheorem{conjecture}[theorem]{Conjecture}%
\usepackage[mode=buildmissing]{standalone}%
    \title{The Universal Property of Measure-Theoretic Probability}\author{
        Eigil Fjeldgren Rischel
      }
\begin{document}
\maketitle

      \begin{abstract}
        Building on work of Chen, we give a universal property of the Markov category BorelStoch of standard Borel spaces and Markov kernels between them. To do this, we introduce a new notion of *coinflip*, or unbiased binary choice, in a Markov category. These are unique if they exist, and automatically preserved by all Markov functors which preserve coproducts. We also provide universal characterizations of various Markov categories of discrete kernels.
      \end{abstract}
    \tableofcontents
    \section{Introduction}

    Markov categories are an abstract approach to probability theory, which axiomatize the structure of categories of probability kernels, beyond the specific details of measure theory. They allow for abstract versions of a number of theorems of classical probability theory, from the 0-1 law of Kolmogorov \cite{rischel-fritz-infinite-products}, the de Finetti theorem \cite{fritz-gonda-perrone-2021}, and a version of the law of large numbers \cite{law-large-nums-fritz-etal}. They also provide an intuitive but rigorous graphical syntax for reasoning about probability.

    In the same way that free Cartesian closed categories describe a minimal syntax for programming with higher-order functions (the simply typed lambda calculus), one may hope that free Markov categories of some sort may provide a minimal syntax for probabilistic programming. In the absence of generating morphisms, the free Markov category on a set of objects is identical with the free Cartesian category on the same set. Hence, although we can consider the free Markov category generated by some collection of primitive distributions, and view its morphisms as terms in a simple probabilistic programming language (Fritz and Liang \cite{fritz-liang-freemarkov-2023} have described the morphisms of this category as certain hypergraphs), the description of these terms as a free Markov category tells us nothing interesting about the equations between different probability distributions. This presents the following question: is there some additional structure on Markov categories such that the free such category has as its morphism probability kernels in the ordinary sense? In this case, the axioms of this structure would be suitable primitives for "simple" probabilistic programming.

    It has been a subject of much speculation whether the "canonical" Markov categories can be described in terms of category theory directly, as a universal object of some kind. Fritz \cite{fritz-stochmat-2009} has given a presentation of \(\mathsf {FinStoch}\) (the Markov category of finite sets and stochastic matrices) by generators and relations, but this simply encodes the set of real numbers into the set of generators. Much more recently, Lorenzin and Zanasi \cite{lorenzin-zanasi-infinite-tensor-2025} have given a description of how to freely adjoin infinite tensor products to a Markov category, and shown that this presents a class of locally constant kernels. Chen \cite{chen-universal-stdborel-2019} has given a universal characterization of the Cartesian category of standard Borel spaces---thus leaving the question of characterizing the class of Markov kernels on top of this.

    Chen's theorem is that the category of standard Borel spaces is initial among extensive, Boolean, countably complete categories (and functors which preserve finite coproducts and countable limits). Clearly \(\mathsf {BorelStoch}\) (or any other interesting Markov category) does not have countable limits, since it lacks products. The right notion of "countable completeness" for Markov categories is given by replacing the infinite products with \emph{Kolmogorov products}, a notion introduced in \cite{rischel-fritz-infinite-products}. Similarly, the "extensive" and "Boolean" parts have to be modified (since we only expect \emph{deterministic} maps into a coproduct to factor over a coproduct decomposition of the domain).

    Beyond these modifications, as noted above, we need some extra structure to provide any nondeterministic morphisms at all. Perhaps the simplest such axiom would be to require the existence of a morphism \(I \to  I + I\), corresponding to a fair coinflip, subject to some equations describing its fairness. It turns out that such a morphism is necessarily unique if it exists, and thus this really becomes a property of, not additional structure on, a Markov category. This property, which we call a \emph{coinflip Markov category}, seems to be about the minimal axiom you could expect to produce a reasonable probability theory.

    With these ingredients, we can state our first theorem:

\begin{theorem}\label{efr-TVLB}

\(\mathsf {BorelStoch}\) is (bi-)initial among countably extensive, Boolean, coinflip Markov categories with countable Kolmogorov products. (And Markov functors which preserve countable coproducts, countable Kolmogorov products and pullbacks along deterministic monomorphisms).

\end{theorem}

    A key inspiration for this theorem is \cite{escardo-simpson-universal-interval-2001}, where Escardo and Simpson give a universal characterization of the interval \([0,1]\) in terms of a "midpoint" operator which is analogous to our binary choice operator. To generate the non-dyadic real numbers they impose beyond basic algebraic properties of the midpoint operator an \emph{iterability} axiom (see \cref{efr-JP37}), which ensures the existence of certain infinite sums. It turns out that in the presence of countable Kolmogorov products this assumption is unnecessary—however, using the iterability axiom, we can also give universal properties to the "discrete" Markov categories like \(\mathsf {FinStoch}\) (which lack these products).

\begin{theorem}\label{efr-FTTL}

  Let \(\kappa \) be a regular cardinal. Let \(\mathsf {Set}_{\bar {\Delta }}^{< \kappa }\) denote the subcategory of \(\mathsf {Set}_{\bar {\Delta }}\) consisting of the sets of cardinality less than \(\kappa \). Then \(\mathsf {Set}_{\bar {\Delta }}^{< \kappa }\) is the initial \(\kappa \)-distributive, \hyperref[efr-Z9IU]{externally iterable}, \hyperref[efr-VTCS]{coinflip} Markov category. In particular, \(\mathsf {FinStoch}\) is the initial distributive, externally iterable, coinflip Markov category, and \(\mathsf {Set}_{\bar {\Delta }}\) is the initial small-distributive, externally iterable, coinflip Markov category.

\end{theorem}

    (Here \(\bar {\Delta }\) denotes the countably supported probability distribution monad). In particular, \(\mathsf {FinStoch}\) is the initial (finitely) distributive, externally iterable, coinflip Markov category, and \(\mathsf {Set}_{\bar {\Delta }}\) is the initial small-distributive, externally iterable, coinflip Markov category.
  
\section{Preliminaries}\label{efr-8P2U}

  The terminology \emph{Markov category} was introduced by Fritz in \cite{fritz-synthetic-markov-cats}, although the idea goes back to the work of Golubtsov, \cite{golubtsov-kleisli}. They have been developed extensively as an abstract foundation for probability theory, see eg \cite{fritz-gonda-perrone-rischel-rep}, \cite{rischel-fritz-infinite-products}, \cite{law-large-nums-fritz-etal}, \cite{fritz-gonda-perrone-2021}. We now review the basics.

\begin{definition}[{Markov category}]

    A \emph{Markov category} \(\mathcal {C}\) is a symmetric monoidal category in which every object is equipped with a comonoid structure \((X, \mathrm {copy}_X, \mathrm {del}_X)\), which is compatible with the monoidal structure in the sense that
    \begin{equation}
      \mathrm {del}_{X \otimes  Y} = \mathrm {del}_X \otimes  \mathrm {del}_Y, \mathrm {copy}_{X \otimes  Y} = (X \otimes  \sigma _{X,Y} \otimes  Y) \circ  (\mathrm {copy}_X \otimes  \mathrm {copy}_Y)
    \end{equation}
    and so that \(\mathrm {del}_X\) is a natural transformation.
  
\end{definition}

  There is a strictification result (\cite{fritz-synthetic-markov-cats}, Theorem 10.17) for Markov categories which allows us to freely elide the structural isomorphisms—we generally do so throughout this paper.

  The morphisms of a Markov category are thought of as probability kernels, that is functions with values in probability measures (where these terms may be interpreted in a suitably abstract sense). The copy and delete maps are the non-stochastic maps which either copy or delete their input.

  It is simple to verify that every Markov category is semiCartesian, and thus the only choice involved in equipping a symmetric monoidal category with a Markov structure is \(\mathrm {copy}_X\). A morphism \(f :X \to  Y \in  \mathcal {C}\) is called \emph{deterministic} if it is a homomorphism for the given comonoid structure. The subcategory of deterministic morphisms is denoted \(\mathcal {C}_\mathrm {det} \subseteq  \mathcal {C}\). It is always a symmetric monoidal subcategory, and Cartesian. Conversely, every Cartesian category carries a unique Markov structure.

  We will often use the notation of Cartesian categories and denote by \(\pi _A\) the map \((A \otimes  \mathrm {del}_B): A \otimes  B \to  A\).

  Given two maps \(f: X \to  A, g: X \to  B\), there is a distinguished pairing \(X \to  A \otimes  B\), given by \((f \otimes  g)\mathrm {copy}_X\). This is called the \emph{independent pairing} of \(f\) and \(g\). Note that by naturality of \(\mathrm {del}\) and the comonoid equations, \(f\) and \(g\) can be recovered from this pairing by postcomposing with \(\pi _A\) or \(\pi _B\). Thus "being independent" is a property of kernels with a tensor product as its codomain. We generalize this to higher tensor products in the obvious way.

  The following class of examples captures almost all Markov categories of interest.

\begin{example}

    If \(\mathcal {C}\) is Cartesian and \(P\) is a monoidal monad on \(\mathcal {C}\), then the Kleisli category (which we denote \(\mathcal {C}_P\)) is a Markov category, with \(\mathrm {copy}_X\) being the diagonal in \(\mathcal {C}\).
  
\end{example}

\begin{definition}

    Let \(\mathcal {C},\mathcal {D}\) be Markov categories. An oplax Markov functor \(F: \mathcal {C} \to  \mathcal {D}\) is an oplax symmetric monoidal functor so that the following diagram commutes:
  
\begin{center}
    \begin {tikzcd}
    F(X) \ar [r] \ar [rd] & F(X \otimes  X) \ar [d] \\
    & F(X) \otimes  F(X)
    \end {tikzcd}
  \end{center}

    An oplax Markov functor is strong if the underlying monoidal functor is strong.
  
\end{definition}

  The undecorated term \emph{Markov functor} was defined in \cite{fritz-synthetic-markov-cats} to mean what we above called a \emph{strong} Markov functor. From this point, "Markov functor" will always mean the strong notion.

  Our goal in this paper is to construct certain functors between Markov categories. The following lemma is helpful in this regard.

\begin{lemma}

    Let \(\mathcal {C}, \mathcal {D}\) be Markov categories and let \(F: \mathcal {C} \to  \mathcal {D}\) be any functor. Then:
  
\begin{enumerate}\item{}\(F\) admits at most one oplax Markov structure, given by the independent pairing of the two projections \(F(\pi _X), F(\pi _Y): F(X \otimes  Y) \to  F(X), F(Y)\).
    \item{}This is a genuine oplax Markov structure if and only if \(F\) preserves deterministic and independent maps, in the sense that if \(X \to  A \otimes  B\) is independent, then \(F(X) \to  F(A \otimes  B) \to  F(A) \otimes  F(B)\) is independent.
    \item{}This oplax Markov structure is strong if and only if the restricted functor \(F_\mathrm {det} : \mathcal {C}_\mathrm {det} \to  \mathcal {D}_\mathrm {det}\) preserves finite products.\end{enumerate}\begin{proof}

      It's straightforward to check that if \(F\) admits any oplax Markov structure, then it preserves deterministic maps and this restricts to an oplax monoidal structure on \(F_\mathrm {det}\) (the fact that the oplaxator itself must be deterministic follows from the associativity). But functors between Cartesian categories always admit a unique oplax monoidal structure. Hence any two oplax structures on \(F\) agree. The only obstruction to this oplax functor extending to a full Markov structure is that it may not be natural with respect to nondeterministic morphisms. Given \(X \to  X', Y \to  Y'\), the tensor \(X \otimes  Y \to  X' \otimes  Y'\) is equal to the independent pairing of \(X \otimes  Y \to  X \to  X'\) and the analogous map to \(Y'\). Applying the assumption that \(F\) preserves independent pairings to this gives naturality. The converse is clear. Finally \(F\) is strong if \(F(X \otimes  Y) \to  F(X) \otimes  F(Y)\) is an isomorphism, but this is exactly the claim that \(F_\mathrm {det}\) preserves products.
    \qedhere
\end{proof}\end{lemma}

    We denote by \(\mathsf {BorelStoch}\) the category of standard Borel spaces and Markov kernels (\cite{fritz-synthetic-markov-cats}, Section 4). We denote by \(\mathsf {Borel}\) the category of standard Borel spaces and measurable maps. Note that \(\mathsf {Borel} = \mathsf {BorelStoch}_\mathrm {det}\). Recall that the standard Borel spaces are those measurable spaces arising as the Borel \(\sigma \)-algebra on Polish spaces. By Kuratowski's theorem these are either discrete on a countable set, or isomorphic to the real numbers (in the Borel \(\sigma \)-algebra). Recall also that \(\mathsf {BorelStoch}\) is the Kleisli category of a monad on \(\mathsf {Borel}\) called the \emph{Giry monad}, see \cite{giry-1982} (see also \cite{fritz-synthetic-markov-cats} section 4 for a review of the history of this concept).

  The categorical notion of product extends obviously to define infinite products. Infinite products exist in many Markov categories of interest, such as those of measurable spaces and Markov kernels, but they clearly can't be characterized by their usual universal property. Instead, we can ask them to be \emph{Kolmogorov} products, in the following sense:

\begin{definition}[{Infinite tensor products}]\label{efr-SIKM}

  Let \(\mathcal {C}\) be a semiCartesian symmetric monoidal category and let \(\{X_j\}_{j \in  J}\) be a family of objects.
  Then if \(P_f(J)\) denotes the set of finite subsets of \(J\) (ordered by inclusion), we obtain a diagram \(P_f(J)^\mathrm {op} \to  \mathcal {C}\), where \(F \mapsto  \bigotimes _{j \in  F} X_j\), and the inclusion \(F \to  F'\) is mapped to the morphism which applies the deletion \(X_j \to  I\) for every \(j\) not in \(F\).

  An \emph{infinite tensor product} of the collection \(\{X_j\}\) is a limit of the cofiltered diagram \(F \mapsto  \bigotimes _{j \in  F} X_j\), if it exists and is preserved by the tensor \(Y \otimes  -\) for every \(Y \in  \mathcal {C}\).

\end{definition}

\begin{definition}[{Kolmogorov product}]\label{efr-EB3U}

  Let \(\mathcal {C}\) be a Markov category. An \hyperref[efr-SIKM]{infinite tensor product} is called a \emph{Kolmogorov product} if all the finite marginals \(\otimes _{j \in  J} X_j \to  \otimes _{j \in  F} X_j\) are deterministic.

\end{definition}

  These were introduced in \cite{rischel-fritz-infinite-products}, where it was also shown that \(\mathsf {BorelStoch}\) admits countable Kolmogorov products (Example 3.6).

\begin{lemma}\label{efr-OH6B}

  Let \(\mathcal {C}\) be a Markov category. Then \(\mathcal {C}\) admits Kolmogorov products of a given cardinality \(\kappa \) if and only if \(\mathcal {C}_\mathrm {det}\) admits Cartesian products of the same cardinality, and the cofiltered limits of \cref{efr-SIKM} are preserved by the inclusion \(\mathcal {C}_\mathrm {det} \to  \mathcal {C}\)
\end{lemma}

  Note that given any Markov functor \(F: \mathcal {C} \to  \mathcal {D}\), and a family \(X_j\) so that the Kolmogorov product \(\bigotimes _j X_j\) exists, there is a family of (necessarily deterministic) maps \(F(\bigotimes _j X_j) \to  \bigotimes _{j \in  F} X_j\) for every finite \(F \subset  J\). As usual we say \(F\) \emph{preserves Kolmogorov products} if this collection exhibits \(F(\bigotimes _j X_j)\) as the Kolmogorov product of the family \(F(X_j)\). If \(\mathcal {C}, \mathcal {D}\) admit Kolmogorov products of a given cardinality, \(F\) preserves them if and only if \(F_\mathrm {det} : \mathcal {C}_\mathrm {det} \to  \mathcal {D}_\mathrm {det}\) preserves products of this cardinality.

  When the family \(X\) is constant, we will denote the Kolmogorov product \(\bigotimes _{i \in  J} X\) simply by \(X^J\).

  If \(\bigotimes _{i \in  J} X_i\) is a Kolmogorov product and \(f_i: A \to  X_i\) is a collection of kernels, there is a canonical independent coupling \((f_i): A \to  \bigotimes _i X_i\), given as the limit of the finite independent couplings \(A \to  \bigotimes _{i \in  F} X_i\). In the case where all \(f_i\) (and all the \(X_i\)) are equal, we will denote this simply as \(f^J : A \to  X^J\).

  The Kolmogorov product \((I + I)^\omega  = 2^\omega \) will play an important role in this paper. We generally write \(0,1 : I \to  2\) for the two canonical points in \(2\). If \(b \in  2^\omega \) is an element, \(0.b, 1.b, 001101.b \in  2^\omega \) denote the result of prepending the given string to \(b\). We may also write \(0.2^\omega \) instead of \(\{0\} \otimes  2^\omega \) for the subobject given by those streams starting with a \(0\), and so on.

\section{Distributive and extensive Markov Categories}\label{efr-MGSE}

  We will need to impose a few different categorical properties on our Markov categories. These are analogues of preexisting properties of ordinary categories, but will typically need to be modified somewhat to suit Markov categories. The usual pattern is to demand that \(\mathcal {C}_\mathrm {det}\) has some property, and that this is preserved by the inclusion into \(\mathcal {C}\). In this section we will describe these conditions.

  For background on extensive and distributive categories, see \cite{carboni-lack-walters-extensive}. For background on Boolean categories, see \cite{johnstone-elephant-vol1}, section A1.4. \cite{chen-universal-stdborel-2019} also contains a review of these terms and the relations between them that suffices for this paper.

  The list of hypotheses we will chiefly be interested in is the following:
  \begin{enumerate}\item{}\(\mathcal {C}_\mathrm {det}\) is Boolean, countably complete and (countably) extensive (so that it receives a unique functor from \(\mathsf {Borel}\) by Chen's theorem).
    
    \item{}
      The inclusion \(\mathcal {C}_\mathrm {det} \hookrightarrow  \mathcal {C}\) preserves the countable coproducts, the pullbacks along coproduct inclusions (hence all monomorphisms), and carries the countable products to Kolmogorov products.
    \end{enumerate}

  There are various ways we can break this up into sub-assumptions. For example, if \(\mathcal {C}\) admits Kolmogorov products of a given cardinality, \(\mathcal {C}_\mathrm {det}\) admits products of that cardinality. So if \(\mathcal {C}\) has countable Kolmogorov products and \(\mathcal {C}_\mathrm {det}\) has pullbacks along monomorphisms, \(\mathcal {C}_\mathrm {det}\) has all limits.

\begin{remark}
\begin{enumerate}\item{}
    In the Cartesian case, if \(\mathcal {C}\) is Boolean, countably complete and extensive, then it is also countably extensive (\cite{chen-universal-stdborel-2019}, Lemma 2.5). However, the proof of this relies on identifying morphisms \(X \to  Y\) with their graphs (subobjects of \(X \times  Y\)), and applying our knowledge about the subobject lattice. In the Markov context, we can't necessarily identify a morphism with a deterministic subobject of \(X\), hence this argument does not work—to be precise, \(\mathcal {C}_\mathrm {det}\) will have all countable coproducts but the inclusion will not necessarily preserve them. If \(\mathcal {C}\) is representable, the inclusion preserves all colimits and this limitation disappears.
  
  \item{}
    Let \(\mathcal {C}, \mathcal {D}\) be Markov categories with countable Kolmogorov products, so that \(\mathcal {C}_\mathrm {det}, \mathcal {D}_\mathrm {det}\) is countably extensive and Boolean, and the inclusions preserve the countable coproducts, and let \(F: \mathcal {C} \to  \mathcal {D}\) be a strong Markov functor. Suppose \(F_\mathrm {det} : \mathcal {C}_\mathrm {det} \to  \mathcal {D}_\mathrm {det}\) preserves countable limits and finite coproducts. Then it automatically preserves countable coproducts as well (again, by Chen's argument). Since these are also countable coproducts in \(\mathcal {C}, \mathcal {D}\), it follows that \(F\) preserves countable coproducts as well.
  \end{enumerate}\end{remark}

\begin{definition}\label{efr-T0AU}

  A Markov category is \emph{distributive} (\(\kappa \)-distributive) if \(\mathcal {C}_\mathrm {det}\) admits finite (\(\kappa \)-small) coproducts, these are preserved by the inclusion \(\mathcal {C}_\mathrm {det} \hookrightarrow  \mathcal {C}\), and by the tensors \(A \otimes  -\).

\end{definition}

  In a distributive monoidal category, if \(A,B\) are each comonoids, then \(A+B\) acquires an induced comonoid structure given by
  \begin{equation}A + B \to  A \otimes  A + B \otimes  B \hookrightarrow  (A+B)\otimes  (A+B)\end{equation}
  If \(\mathcal {C}\) is a Markov category, each of these objects has furthermore a \emph{canonical} given comonoid structure. \(\mathcal {C}\) is distributive \emph{as a Markov category} if and only if the coproducts can be chosen so that the canonical comonoid structure coincides with the induced one.

  This parametrization in a regular cardinal \(\kappa \) will recur a few times in this paper.
  We are chiefly interested in the cases \(\kappa  = \omega \), corresponding to \emph{finite} coproducts, and \(\kappa  = \omega _1\), corresponding to \emph{countable} coproducts (this case requires a weak form of the axiom of choice).

\begin{definition}[{Extensive Markov Category}]\label{efr-B2P7}

  A Markov category is said to be an \emph{extensive Markov category} if it admits finite coproducts, whose injections are deterministic, and which satisfy the following equivalent conditions:

\begin{enumerate}\item{}If we let \(\mathcal {C}_{/a}^\mathrm {det}\) refer to the full subcategory of the slice spanned by the deterministic morphism \(x \to  a\), we have an equivalence of categories \(\mathcal {C}_{/a}^\mathrm {det} \times  \mathcal {C}_{/b}^\mathrm {det} \cong  \mathcal {C}_{/a + b}^\mathrm {det}\), given by taking coproducts
  \item{}\(\mathcal {C}_\mathrm {det}\) is an extensive category in the usual sense and the inclusion \(\mathcal {C}_\mathrm {det} \to  \mathcal {C}\) preserves pullbacks along coproduct inclusions.\end{enumerate}

  More generally, for an infinite regular cardinal \(\kappa \), we say that \(\mathcal {C}\) is \(\kappa \)-extensive if \(\mathcal {C}_\mathrm {det}\) is a \(\kappa \)-extensive category in the ordinary sense and both coproduct inclusions and pullbacks along them are preserved by the functor \(\mathcal {C}_\mathrm {det} \to  \mathcal {C}\).

\end{definition}

  Note that any extensive Markov category is automatically distributive (for the same \(\kappa \)). Note also that if \(\mathcal {C}\) is extensive Markov, then \(\mathcal {C}_\mathrm {det}\) automatically admits all finite limits, since it has finite products just by virtue of being a Markov category.

  Recall that a coherent category is called \emph{Boolean} if every subobject lattice is a Boolean algebra. Also note that an extensive category is Boolean (and coherent) as soon as it satisfies the condition that every monomorphism \(V \hookrightarrow  X\) is a coproduct inclusion (the object \(V'\) so that \(A = V + V'\) is then the complement subobject of \(V\)). We will not need to study any other case than this, and so we simply make the following definition.

\begin{definition}[{Boolean Markov category}]\label{efr-7JNJ}

  An extensive Markov category \(\mathcal {C}\) is \emph{Boolean} (as a Markov category) if \(\mathcal {C}_\mathrm {det}\) is Boolean, i.e if every monomorphism in \(\mathcal {C}_\mathrm {det}\) is a coproduct inclusion.

\end{definition}

  In the common case that \(\mathcal {C}\) is the Kleisli category of a monad on \(\mathcal {C}_\mathrm {det}\), the inclusion \(\mathcal {C}_\mathrm {det} \hookrightarrow  \mathcal {C}\) automatically preserves all colimits, which makes many of the above things automatic.

\begin{lemma}\label{efr-198L}

    Let \(\mathcal {C}\) be a \emph{representable} Markov category, and let \(P\) be the associated distribution monad. Then

\begin{enumerate}\item{}
      If \(\mathcal {C}_\mathrm {det}\) is \(\kappa \)-extensive and \(P: \mathcal {C}_\mathrm {det} \to  \mathcal {C}_\mathrm {det}\) preserves pullbacks along coproduct inclusions, then \(\mathcal {C}\) is \(\kappa \)-extensive.
    
    \item{}
      If \(\mathcal {C}_\mathrm {det}\) is Boolean, then \(\mathcal {C}\) is Boolean.
    
    \item{}
      If \(\mathcal {C}_\mathrm {det}\) is (finitely) extensive, Boolean, and \(\kappa \)-complete, and \(P: \mathcal {C}_\mathrm {det} \to  \mathcal {C}_\mathrm {det}\) preserves the cofiltered limits \(\prod _{i \in  J} A_i = \lim _{F \subseteq  J \mathrm { finite}} \prod _{i \in  F} A_i\) for \(|J| < \kappa \), then \(\mathcal {C}\) is \(\kappa \)-extensive (and Boolean and admits \(\kappa \)-small Kolmogorov products).
    \end{enumerate}\end{lemma}

\begin{remark}\label{efr-GZCD}

  We will characterize \(\mathsf {BorelStoch}\) as the initial Markov category which is countably extensive, Boolean, admits countable Kolmogorov products, and a \hyperref[efr-0XH9]{coinflip}. Let us say a few words justifying these axioms as natural. Extensivity is a property of most categories of "spaces", and the requirement that \(\mathcal {C}_\mathrm {det} \hookrightarrow  \mathcal {C}\) preserves coproducts and pullbacks along monomorphisms holds in essentially every Markov category of interest.

  However, Boolean categories are a bit more uncommon. In the presence of extensivity, being Boolean is equivalent to \(I + I\) being a subobject classifier - the obstruction to this typically being that "the" indicator of many subobjects are discontinuous.

  A common assumption about Markov categories is the presence of \emph{conditionals} (\cite{fritz-synthetic-markov-cats} 11.5). If \(V \hookrightarrow  X\) is a subobject, and \(I \to  I + I, I \to  V, I \to  X\) are distributions of full support, then a Bayesian inverse of their pairing \(I + I \to  X\) is forced to have a discontinuity at the boundary of \(V\).

  While this does not rise to the level of a formal argument that conditionals imply Boolean (it seems very hard to imagine such an argument, since conditionals are only defined up to almost sure equality, and measures with full support don't always exist, for example, in \(\mathsf {BorelStoch}\)), it does demonstrate the difficulties involved in admitting conditionals without being Boolean.

\end{remark}

\begin{proposition}\label{efr-GY4W}

\(\mathsf {BorelStoch}\) is countably extensive and Boolean.

\begin{proof}

    Note first that \(\mathsf {Borel}\) is countably extensive and Boolean, by \cite{chen-universal-stdborel-2019} theorem 1.1. Hence it suffices to observe that the Giry monad preserves pullbacks along coproduct inclusions. Consider \(X \to  A + B\). The claim is that the square
  
\begin{center}
    \begin {tikzcd}
    G(X_A) \ar [r] \ar [d] & G(A) \ar [d] \\
    G(X) \ar [r] & G(A + B)
    \end {tikzcd}
  \end{center}

    is a pullback, where \(X_A\) is the preimage of \(A\) inside \(X\). But this is clear: a probability measure on \(A\) is exactly a probability measure on \(A + B\) which happens to be concentrated on \(A\), and a probability measure on \(X\) has an image in \(G(A + B)\) of this form if and only if it is concentrated on \(X_A\) (and is such equivalent to a measure on \(X_A\)).
  \qedhere
\end{proof}\end{proposition}
\section{Midpoint algebras and Coinflips}\label{efr-AHFX}

    We will not have any hope of giving a universal property for \(\mathsf {BorelStoch}\), or any interesting Markov category, without some property that forces certain maps to be nondeterministic. We now give such a structure. As mentioned, our basic idea will be to impose the existence of a "unbiased binary random choice". In a general Markov category, this choice forms a morphism \(A \otimes  A \to  A\). The axioms that such a choice must satisfy are essentially the axioms described by Escardo and Simpson in \cite{escardo-simpson-universal-interval-2001}, which they termed \emph{midpoint algebras}. The main theorem of their paper is a characterization of the interval \([0,1]\) as the free iterable (\cref{efr-JP37}) midpoint algebra on two points—this will play a central role in this paper.

\begin{definition}[{Midpoint algebra}]\label{efr-6W7X}

  Let \(\mathcal {C}\) be a Cartesian category. Then a \emph{midpoint algebra} is an object \(A\) equipped with \(m: A \times  A \to  A\), so that:
  \begin{equation}m(a,a) = a\end{equation}
  \begin{equation}m(a,b) = m(b,a)\end{equation}
  \begin{equation}m(m(a,b),m(c,d)) = m(m(a,c), m(b,d))\end{equation}

  We will call a homomorphism of midpoint algebras (a morphism satisfying \(m(f(a),f(a')) = f(m(a,a'))\)) a \emph{midpoint homomorphism}
\end{definition}

\begin{definition}[{Internal midpoint algebra}]\label{efr-1XB5}

  Let \(\mathcal {C}\) be any category. An \emph{internal midpoint algebra} structure on an object \(A \in  \mathcal {C}\) is a family of midpoint algebra structures on each \(\mathcal {C}(X,A)\), so that precomposition is a midpoint homomorphism for every \(f: Y \to  X\).

\end{definition}

    Note that if \(\mathcal {C}\) is Cartesian this agrees with the above definition.

\begin{definition}[{Coinflip structure}]\label{efr-QBVU}

  Let \(\mathcal {C}\) be a Markov category. A \emph{coinflip structure} is a \hyperref[efr-1XB5]{internal midpoint algebra} structure on each \(A \in  \mathcal {C}\), which is natural in \(A\). In other words, it is a \hyperref[efr-6W7X]{midpoint algebra} structure on each homset \(\mathcal {C}(X,Y)\) so that both pre- and postcomposition are midpoint homomorphisms.

\end{definition}

\begin{proposition}\label{efr-OX26}

  Let \(\mathcal {C}\) be a Markov category. Then every coinflip structure on \(\mathcal {C}\) is equal.

\begin{proof}

    Let \(m_1,m_2\) denote two coinflip structures. Let \(\mu _1 = m_1(\pi _1, \pi _2) : A \otimes  A \to  A\), analogously \(\mu _2\) (we condense the notation by writing \(\mu _1,\mu _2\) regardless of the object in question). Note \(m_1(a,b) = \mu _1\langle  a,b \rangle \), where \(\langle  a,b \rangle \) is \emph{any} pairing of \(a,b\), by naturality.

    First consider \(m_1(m_2(a,b),m_2(c,d))\). By the above this is \(\mu _1 \langle  \mu _2 \langle  a,b \rangle , \mu _2 \langle  c,d \rangle  \rangle \). Applying naturality to the morphism \(\mu _2\), we can also rewrite this as \(\mu _2 \circ  (m_1(\langle  a,c \rangle ,\langle  b,d \rangle ))\)

    Now take \(c = b, d = a\). By symmetry for \(m_2\) we can rewrite the left-hand side as \(m_1(m_2(a,b),m_2(a,b))\). By idempotency this is simply \(m_2(a,b)\). The other expression now becomes \(\mu _2 \circ  (m_1(\langle  a,b \rangle ,\langle  b,a \rangle ))\). Now consider the first projection of the map \(m_1(\langle  a,b \rangle ,\langle  b,a \rangle )\). By naturality it is equal to \(m_1(a,b)\). By symmetry, so is the second projection. Hence this is a pairing of \(m_1(a,b)\) with itself. Hence the composite is equal to \(m_2(m_1(a,b),m_1(a,b)) = m_1(a,b)\). This concludes the proof.
  \qedhere
\end{proof}\end{proposition}

\begin{definition}[{Coinflip Markov Category}]\label{efr-VTCS}

  A \emph{coinflip Markov category} is a Markov category which admits a coinflip structure.

\end{definition}

    Note that any coinflip structure is determined by the family of maps \(m(\pi _0,\pi _1) : A \otimes  A \to  A\). To define a coinflip structure, such a family \(\mu _A\) must satisfy the equations of a midpoint algebra, be natural in \(A\), and have the further property that \(\mu _A f = \mu _A (\pi _0 f, \pi _1 f)\) (that is, the composite depends only on the projections of \(f\)).

    Coinflip structures are not preserved by all Markov functors, since \(F(\mu _A)\) only has to enjoy the above property with respect to maps in the image of \(F\). However, as we now show, this is true in the distributive case.

\begin{definition}[{Coinflip in a Markov category}]\label{efr-0XH9}

  Let \(\mathcal {C}\) be a Markov category with finite coproducts. A \emph{coinflip} in \(\mathcal {C}\) is a morphism \(f: I \to  I + I\) obeying the following equations:
  \begin{enumerate}\item{}
      Symmetry: \(\tau _{I,I} f = f\), where \(\tau _{A,B} :A + B \to  B + A\) is the canonical isomorphism. 
    
    \item{}
      Interchange: \((f + f)f = (I + \tau _{I,I} + I)(f + f)f\)\end{enumerate}

  Given a coinflip in a distributive Markov category, there is a canonical induced natural transformation \(A \otimes  A \to  A\), which makes every \(A \in  \mathcal {C}\) into a midpoint algebra.

\end{definition}

\begin{proposition}\label{efr-V3JG}

  Let \(\mathcal {C}\) be a distributive Markov category. Any two coinflips in \(\mathcal {C}\) coincide. In particular, any Markov functor \(\mathcal {C} \to  \mathcal {D}\) which preserves these finite coproducts also preserves the coinflip.

\begin{proof}

    Let \(f\) be a coinflip. Consider the map \(2 \xrightarrow {\oplus  f} 2\) given by flipping the coin, applying \(\neg : 2 \to  2\) in one branch and doing nothing in the other. We claim this is equal to deletion followed by the coinflip. This can be checked on each summand (since \(2\) is a coproduct), where it follows by symmetry.

    Now let \(f_1, f_2\) be two coinflips and consider the operation given by \(1 \xrightarrow {f_1 \otimes  f_2} 2 \otimes  2 \xrightarrow {\oplus }\). This is easily seen to be equal both to \(f_1 ; (\oplus  f_2)\) and \(f_2 ; (\oplus  f_1)\). But these are equal to \(f_2\) and \(f_1\) respectively by the above.
  \qedhere
\end{proof}\end{proposition}

    Note that any coinflip \(I \to  I + I\) in a distributive Markov category induces a coinflip structure, and vice versa. Since both notions are unique, this is obviously a 1-1 correspondence. From this we easily derive the following:

\begin{corollary}\label{efr-7QV9}

  Let \(\mathcal {C}, \mathcal {D}\) be distributive coinflip Markov categories, and suppose \(F: \mathcal {C} \to  \mathcal {D}\) is a Markov functor which preserves coproducts and the monoidal unit. Then \(F\) also preserves the coinflip structure.

\end{corollary}

    The terminology "midpoint algebra" is from \cite{escardo-simpson-universal-interval-2001}. In that paper, Escardo and Simpson provide a characterization of the interval as the free midpoint algebra satisfying certain properties generated by two points. The key property is the following:

    Escardo and Simpson introduced the following notion:

\begin{definition}[{Iterable midpoint algebra}]\label{efr-JP37}

  A \hyperref[efr-6W7X]{midpoint algebra} \((A,m)\) in \(\mathcal {C}\) is \emph{iterable} if, for each \(s: X \to  X \times  A \in  \mathcal {C}\), there is a unique \(u : X \to  A\) so that the diagram
  
  \begin{center}
    \begin {tikzcd}
    X \ar [d, "u"] \ar [r, "s"] & X \times  A \ar [d, "u \times  A"] \\
    A & A \times  A \ar [l, "m"]
    \end {tikzcd}
  \end{center}
\end{definition}

    If \(A\) is an internal midpoint algebra in a category \(\mathcal {C}\), a priori, the only version of \cref{efr-JP37} we can impose is that all the midpoint algebras \(\mathcal {C}(X,A)\) are iterable (in \(\mathsf {Set}\)). However, if \(\mathcal {C}\) is Markov, we can alternatively demand the following (which is just \cref{efr-JP37} with tensor products instead of products):

\begin{definition}[{Internally iterable}]\label{efr-87W1}

  Let \(\mathcal {C}\) be a Markov category, let \(A\) be an object, and let \(m\) be a midpoint algebra on \(A\), i.e a midpoint algebra structure on each \(\mathcal {C}(X,A)\), natural in \(X\). Then \(A\) is \emph{internally iterable} if, for each \(X \to  X \otimes  A \in  \mathcal {C}\), there exists a unique map \(u: X \to  A\) satisfying \(u = m(u \pi _X, \pi _A)s\)
\end{definition}

    In this case we might call \(u\) an \emph{iteration operator} for \(s\).

\begin{definition}\label{efr-Z9IU}

  Let \(\mathcal {C}\) be a coinflip Markov category. Then it is called \emph{externally iterable} if the midpoint algebras \(\mathcal {C}(X,Y)\) are all iterable in \(\mathsf {Set}\). It is said to be \emph{internally iterable} if they are all internally iterable.

\end{definition}

    Escardo and Simpson's theorem is that the set \([0,1]\) is the free iterable, \emph{cancellative} (i.e satisfying \(m(a,b) = m(a,c) \Rightarrow  b = c\)) midpoint algebra on two generators. We now strengthen their theorem by proving that it is free among merely iterable midpoint algebras (since it is cancellative, this implies their statement).

\begin{lemma}\label{efr-CZUI}

  In \(\mathsf {Set}\), the free iterable midpoint algebra on two generators \(0,1\) is \([0,1]\). That is, given an iterable midpoint algebra \(A\) and two points \(a_0, a_1\), there exists a unique midpoint homomorphism \(f: [0,1] \to  A\) so that \(f(0)=a_0, f(1)=a_1\).

\begin{proof}

    Let \(A\) be an iterable midpoint algebra. Note that any map \(f: [0,1] \to  A\) must be a homomorphism for the infinitary midpoint operator \(M\). Hence given \(a_0,a_1 \in  A\), the unique \(f\) so that \(f(0) = a_0, f(1) = a_1\) is given by \(f(0.b_1b_2 \dots ) = M(a_{b_1}, \dots )\). This is straightforwardly seen to be well-defined and a midpoint homomorphism.
  \qedhere
\end{proof}\end{lemma}

\begin{lemma}\label{efr-JB06}

  Let \(X\) be a possibly infinite set. Then the infinitary simplex \(\bar {\Delta }(X) =  \{ (a_x \geq  0)_{x \in  X} \mid  \sum  a_x = 1 \}\) is the free iterable midpoint algebra on \(X\).

\begin{proof}

    We identify the elements of \(\bar {\Delta }(X)\) with countably supported probability measures on \(X\). Observe that \(\bar {\Delta } X\) is an iterable midpoint algebra (with the infinitary choice operation given by the sum \(M((a_i)) = \sum _i 1/2^i a_i\)). By the preceding, every finitely supported measure can be obtained by iterating \(M\), and every countably supported measure can be obtained by applying \(M\) to a sequence of these (half of the probability measure must be concentrated in some finite subset, then \(1/4\) of the remainder again in some finite subset, and so on). Hence there is at most one midpoint homomorphism \(\bar {\Delta }X \to  Y\) extending any given function \(f: X \to  Y\), for any iterable midpoint algebra \(Y\).

    Now we must show that one exists. By a \emph{finitary dyadic measure}, we mean one which is finitely supported and where each weight \(a_x\) has the form \(k/2^n\). Note that these are exactly those that can be written using just the binary choice operator \(m(x,y)\). Also note that any two such terms are equal in a generic midpoint algebra if and only if they represent the same such measure. Hence given \(f: X \to  Y\), there is a unique and well-defined \(f(a)\) defined on the finitary dyadic measures, which is a midpoint homomorphism.

    Now for any probability measure \(\mu \), we can write it as \(\mu  = M(a_1,\dots )\), where each \(a_i\) is a finitary dyadic measure. We define our operation as \(f(\mu ) = M(f(a_1), \dots )\), with \(f(a_i)\) defined as above. It suffices to show that this is well-defined, since 
    \begin{equation}f(m(\mu ^1,\mu ^2)) = f(m(M(a_1^1,\dots ),M(a_1^2,\dots ))) =\end{equation}
    \begin{equation}f(M(m(a^1_1,a^2_1), \dots )) = M(m(f(a_1^1),f(a_1^2)),\dots ) = m(f(\mu ^1),f(\mu ^2))\end{equation}

    Now let \(M(a_1,\dots ) = M(b_1,\dots ) = \mu \) be two representations in this form of the same measure. Then there exists some smallest \(N\) so that there exists a finitary dyadic probability measure \(c_1 = \sum _i \frac {k_i}{2^N} \delta _{x_i}\) so that \(k_i/2^N < 2 \mu (x_i)\) whenever $\mu(x_i)$ is positive - that is, so that \(c/2\) is \emph{strictly} less than \(\mu \) everywhere on the support.

    Now in the sum \(\sum _i 1/2^i a_i\) must exceed \(c/2\) by some \emph{finite} partial sum, say \(M\). Then \(\sum ^M_{i=1} 1/2^i a_i\) is a finitary dyadic measure. Hence a finite manipulation using the properties of midpoint algebras can rewrite the term \(M(a_1,\dots )\) into \(M(c_1,a_2',\dots ,a_M',a_{M+1})\). Iterating this we build \(c_2,\dots \) so that \(M(a_1,\dots ) = M(c_1,\dots )\). Since these only depended on the measure \(\mu \), we similarly find \(M(b_1,\dots ) = M(c_1,\dots )\). Since this argument used only the axioms of an iterable midpoint algebra, it follows that \(f\) is well-defined.
  \qedhere
\end{proof}\end{lemma}

\begin{corollary}\label{efr-REOY}

  Any iterable midpoint algebra carries the structure of an algebra of \(\bar {\Delta }\).
  Moreover any midpoint homomorphism between iterable midpoint algebras is a \(\bar {\Delta }\)-homomorphism. In particular, if \(\mathcal {C}\) is an externally iterable coinflip Markov category, each homset \(\mathcal {C}(X,Y)\) is a \(\bar {\Delta }\)-algebra, and composition is biconvex (i.e convex in each variable separately). If \(\mathcal {C} \to  \mathcal {D}\) is a functor between externally iterable coinflip Markov categories which preserves the coinflip (eg if \(\mathcal {C}, \mathcal {D}\) are distributive), it automatically preserves this structure.

\end{corollary}

    The situation with coinflip Markov categories is similar to the situation of preadditive categories, known for example from the study of homological algebra. Under mild conditions on a category (admitting finite products and coproducts is sufficient), an enrichment over commutative monoids is unique if it exists.

\begin{proposition}\label{efr-E4WU}

  Let \(\mathcal {C}\) be a coinflip Markov category. If it is externally iterable, then it is internally iterable.

\begin{proof}

    Let \(s: X \to  A \otimes  X\) be a map, and suppose \(u_1,u_2\) both satisfy the equation \(u_i = m(\pi _A s, u_i \pi _X s)\). Let \(a_n : X \to  A\) denote the map defined inductively by \(a_1 = \pi _A s\), \(a_{n+1} = a_n \pi _X s\). Let \(b_i^{n}\) be defined inductively by \(b_i^1 = u_i\), \(b_i^{n+1} = b_i^n \pi _X s\). Then we see that \(b_i^n = m(a_n, b_i^{n+1})\) for both \(i=1,2\). Therefore they must agree, and in particular \(u_1 = u_2\). This concludes the proof.
  \qedhere
\end{proof}\end{proposition}

\begin{proposition}\label{efr-60YW}

  Let \(\mathcal {C}\) be a coinflip Markov category which is countably distributive. If \(\mathcal {C}\) is internally iterable, then it is externally iterable.

\begin{proof}

    By \cite{escardo-simpson-universal-interval-2001}, a midpoint algebra \((A,m)\) in \(\mathsf {Set}\) is iterable if and only if is satisfies these conditions:
    \begin{enumerate}\item{}
        There exists an operator \(M: A^\omega  \to  A\) so that \(M(a_1,a_2, \dots ) = m(a_1, M(a_2, a_3, \dots ))\)
      \item{}
        Given two sequences \(a_i, b_i\) so that \(b_i = m(a_i, b_{i+1})\) for all \(i\), then \(b_1 = M(a_1, a_2, \dots )\)\end{enumerate}

    Let a countable family \(f_i : A \to  B\) be given. Consider the operation \(\sum ^\omega _{i=0} A \to  B \times  \sum ^\omega  A\) given by the copairing of all the \(f_i\) paired with the map \(\sum ^\omega _{i=0} A \to  \sum ^\omega _{i=0} A\) which maps each summand to the next one (i.e \((i,a) \mapsto  (i+1,a)\)). Now an iteration of this map \(u: \sum _i A \to  B\) is equivalently a sequence of maps \(g_i\) so that \(g_i = m(f_i,g_{i+1})\). Hence the existence and uniqueness of such an operator precisely corresponds to points 1 and 2 above.
  \qedhere
\end{proof}\end{proposition}

    If the Markov category under consideration has at least countable coproducts, we may therefore simply speak of an \emph{iterable} coinflip Markov category.
  
\section{The universal property of discrete probability}\label{efr-ALKX}

  Using our description of the free iterable midpoint algebras, we can now prove \cref{efr-FTTL}.

\begin{proof}

  Let \(\mathcal {C}\) be a \(\kappa \)-distributive externally iterable coinflip Markov category. There is an essentially unique \(\kappa \)-coproduct preserving Markov functor \(F: \mathsf {Set}^{< \kappa } \to  \mathcal {C}\) (since \(\mathcal {C}_\mathrm {det}\) is \(\kappa \)-distributive and \(\mathsf {Set}^{< \kappa }\) is initial such). Any extension to the stochastic morphisms is determined by its action on the homsets \(\mathsf {Set}_{\bar {\Delta }}^{< \kappa }(*,X) = \bar {\Delta }(X)\), by the coproduct universal property. On these its action must be given by taking a convex combination \(\sum _i \epsilon _i x_i\) to the convex combination \(\sum _i \epsilon _i F(x_i)\). This is functorial by \cref{efr-REOY}, finishing the proof.
\qedhere
\end{proof}

  We can also describe the free iterable coinflip Markov category on a generic category, as follows:

\begin{proposition}

    Let \(\mathcal {C}\) be a Markov category. Let \(\bar {\Delta } (\mathcal {C})\) denote the category with the same objects, whose homsets are given by \(\bar {\Delta } (\mathcal {C}(X,Y))\), with the uniquely extended bilinear composition. Then \(\bar {\Delta } \mathcal {C}\) is a Markov category, and \(\mathcal {C} \hookrightarrow  \bar {\Delta } \mathcal {C}\) presents it as the free externally iterable coinflip Markov category generated by \(\mathcal {C}\)
\end{proposition}

  Note that, although coinflip structures are unique, not every functor preserves them. Hence the free coinflip Markov category on a category which is already coinflip is not the category itself, as one might expect. It is interesting to consider whether a version of this theorem for distributive Markov categories exists (which would generalize \cref{efr-FTTL}), but we do not currently know how to prove such a theorem.

\section{The universal property of \(\mathsf {BorelStoch}\)}\label{efr-D09F}

    We are now ready to prove \cref{efr-TVLB}. The basic point is that any Markov kernel between standard Borel spaces can be written as the composite of a measurable function \(f: A \times  2^\omega  \to  B\) with the uniform measure on \(c^\omega : I \to  2^\omega \). By Chen's theorem, the image of \(f\) is uniquely determined, and by \cref{efr-V3JG}, so is the image of \(c^\omega \). Hence there is at most one such functor---what we need to prove is that the above is well-defined (that is, choosing different factorizations gives the same kernel in the image) and actually forms a Markov functor.

    The strategy for the proof is as follows:
  
\begin{enumerate}\item{}
      Observe that for any \(\mathcal {C}\), there is an essentially unique functor \(\mathsf {Borel} \to  \mathcal {C}\) (which lands in \(\mathcal {C}_\mathrm {det}\)), by Chen's theorem.
    
    \item{}
      Prove that there is a (Markov) extension to the category of discrete kernels, \(\mathsf {Borel}_{\bar {\Delta }}\).
    
    \item{}
      Finally prove that there is a further Markov extension to \(\mathsf {BorelStoch}\)\end{enumerate}

    We have essentially already given the proof of the first step in the introduction, but let us state it here for completeness.

\begin{proposition}\label{efr-AFCK}

  Let \(\mathcal {C}\) be a Markov category which is \hyperref[efr-B2P7]{countably extensive}, \hyperref[efr-7JNJ]{Boolean}, \hyperref[efr-VTCS]{coinflip}, and admits countable \hyperref[efr-EB3U]{Kolmogorov products}. Then there is an essentially unique (strong) Markov functor \(i: \mathsf {Borel} \to  \mathcal {C}\) which preserves pullbacks along monomorphisms and countable Kolmogorov products.

\begin{proof}

    Since \(\mathsf {Borel}\) is Cartesian, any such Markov functor must factor over \(\mathcal {C}_\mathrm {det}\). \(\mathcal {C}_\mathrm {det}\) has countable limits (it has countable products by \cref{efr-OH6B} and pullbacks along monomorphisms, which along with finite products suffices to build equalizers) and is countably extensive and Boolean by assumption. The assumptions on \(i\) are equivalent to the factorization \(\mathsf {Borel} \to  \mathcal {C}_\mathrm {det}\) preserving countable limits and countable coproducts. Hence this uniqueness is exactly \cite{chen-universal-stdborel-2019} theorem 1.1
  \qedhere
\end{proof}\end{proposition}

    To make our results slightly more broadly applicable, since there aren't many Boolean categories, we will instead take a suitable Markov functor \(\mathsf {Borel} \to  \mathcal {C}\) as the input, and prove that given such, there is a unique extension of the above kind. In the below, we will generally abuse notation and identify a standard Borel space with its image under this functor.

    Let's now see the details of part 2. We begin with the following lemma, which is what allows us to use infinitary sampling without the iterability assumption.

\begin{lemma}\label{efr-K4NK}

  Let \(\mathcal {C}\) be an extensive coinflip Markov category with countable Kolmogorov products, and let \(\mathsf {Borel} \to  \mathcal {C}\) be a strong Markov functor which preserves pullbacks along coproduct inclusions and countable Kolmogorov products. Let \(V_1 \subseteq  2^\omega \) be the subobject consisting of those sequences containing infinitely many ones. Then the map \(c^\omega  : I \to  2^\omega \) factors over \(V_1\). In particular, it factors over the inclusion \(1 . 2^\omega  \sqcup  01.2^\omega  \sqcup  \cdots  \hookrightarrow  2^\omega \) of the subobject of sequences with at least one \(1\).

\begin{proof}

    Let \(r_1: 2^\omega  \to  2\) be the indicator of the subset consisting of those streams which have only finitely many \(1\)s. Observe that for any finite \(N\), \(r_1\) factors over the second coordinate of the decomposition \(2^\omega  \to  2^N \times  2^\omega \).

    As a result, if we form the joint \(I \to  2^\omega  \otimes  2\) given by sampling the uniform distribution and computing \(r_1\), the second coordinate is independent of any finite prefix of the first. Hence, by the abstract Kolmogorov zero-one law (\cite{rischel-fritz-infinite-products}, theorem 5.3), the composite \(r_1 c^\omega  : I \to  I+I\) is deterministic.

    Hence it classifies some subterminal object \(q\). By symmetry, \(r_0 c^\omega \) is the indicator of \(q\) as well, where \(r_0\) is the indicator of sequences with only finite many \(0\)s. But these two subterminal objects must be disjoint, hence they must both be empty. (This step uses the fact that pullbacks along monomorphisms in \(\mathcal {C}_\mathrm {det}\) are still pullbacks in \(\mathcal {C}\))

    So \(r_0c^\omega , r_1c^\omega \) are both equal to \(0 : I \to  2\). It follows that \(c^\omega \) factors over the subobject \(2^\omega  \setminus  \{0\}^\omega  \hookrightarrow  2^\omega \) (since any sequence with finitely many zeroes certainly has at least one \(1\)).
  \qedhere
\end{proof}\end{lemma}

    Thus we can define a countable-arity operator \(M_A: A^\omega  \to  A\) in any of these categories by drawing a sample from \(c^\omega  \in  2^\omega \) and projecting onto the coordinate of the first \(1\) (which necessarily exists).

      For any standard Borel space \(A\), let \(BT(A)\) denote the (again standard Borel) space of balanced, finite binary trees with leaves labeled by \(A\). This can be identified with the coproduct \(\sum _{n \in  \mathbb {N}} A^{2^n}\). Let \(g : BT(A) \times  BT(A) \to  BT(A)\) glue two trees together, duplicating the shallower one so that they have the same depth. In \(\mathcal {C}\), there is a unique kernel \(BT(A) \to  A\) which carries this operation to the midpoint operation in \(\mathcal {C}\). Note that this can be factored as a map \(BT(A) \times  2^\omega  \to  A\) which uses the first \(n\) bits to choose a branch, composed with the uniform distribution.

    Given anything like \((BT(A)^{\omega  \times  \omega })^{\omega }\), we can build a well-defined map into \(A\) by composing these sampling homomorphisms. This is what we mean below by a sentence like "the sampling homomorphism \(BT(A)^{\omega ^N} \to  A\). Note that this is well-defined in any extensive, coinflip \(\mathcal {C}\) with countable Kolmogorov products equipped with \(\mathsf {Borel} \to  \mathcal {C}\).

    Note that we have not proven that \(M_A\) is uniquely determined by the coinductive definition \(M_A(a_1, \dots ) = m(a_1, M(a_2, \dots ))\), merely constructed a canonical such operation. We have not been able to prove uniqueness in general (which would imply that any countably distributive, Boolean, countably complete coinflip Markov category is iterable), but neither do we know of a counterexample.

\begin{lemma}\label{efr-X0VA}

  With \(\mathcal {C}, \mathsf {Borel} \to  \mathcal {C}\) as above, let \(A\) be (the image of) a standard Borel space, and consider the map \(\lambda : (BT(A)^\omega )^\omega  \to  A\) given by sampling two natural numbers according to \(c^\omega : I \to  2^\omega  \setminus  0 \to  \omega \) and indexing with them, then using the canonical sampling map. There exists a map \(BT(A)^{\omega  \times  \omega } \to  BT(A)^\omega \) which preserves the sampling map. Moreover this map also preserves the map into \(\Delta (A)\).

\begin{proof}

    Let \(F: BT(A)^{\omega  \times  \omega } \times  (2^\omega  \setminus  0)^2 \times  (2^\omega )^\omega  \to  A\) be a parametrization of the sampling map. Specifically, the first two bitstreams select the indexes \((i,j)\) into the matrix of finitary distributions, after which the \(i\)th of the remaining bitstreams is used to choose a branch.

    Let us adopt the convention that the first index chooses the \emph{column} and the second index chooses the \emph{row} when describing the manipulations below.

    Given a sequence \((d_{ij}) \in  BT(A)^{\omega  \times  \omega }\), construct \(r_1(d)_{ij}\) as follows:
  
\begin{equation}r_1(d)_{1j} = g(d_{11}, g(d_{12}, d_{21}))\end{equation}
\begin{equation}r_1(d)_{2j} = g(d_{1(j+2)},d_{2(j+1)})\end{equation}
\begin{equation}r_1(d)_{ij} = d_ij \text{ for $i \geq 3$}\end{equation}

    First note that clearly \(r_1(d)\) has the same image in \(\Delta (A)\) as \(d\)—the first column of \(r_1(d)\) has the same sum as the three most probably indices of \(d\), and the second column of \(r_1(d)\) has the same sum as the remainder of the first two columns of \(d\). Also note that there exists a permutation \(\rho _1 : (2^\omega  \setminus 0 )^2 \times  (2^\omega )^omega  \to   (2^\omega  \setminus  0)^2 \times  (2^\omega )^omega \) so that \(F(r_1(d),\rho _1(b)) = F(d,b)\).

\(\rho _1(b_1,b_2,(b_3^i))\) can be chosen as follows: it does not modify \(b_3^i\) for \(i>2\) (since these bits are only used for choices in the first column). Omitting therefore these strings from the notation below, we define
    \begin{equation}\rho _1(1.c,1.r,s^1,s^2) = (1.c,r,1.s^1,s^2)\end{equation}
    \begin{equation}\rho _1(1.c,01.r,s^1,s^2) = (1.c, r,01.s^1,s^2)\end{equation}
    \begin{equation}\rho _1(01.c,1.r,s^1,s^2) = (1.c,r,10.s^2, s^1)\end{equation}
    \begin{equation}\rho _1(1.c, 00.r, s^1,s^2) = (01.c, r, s^2, 1.s^1)\end{equation}
    \begin{equation}\rho _1(01.c, 0.r, s^1,s^2) = (01.c, r, s^1, 0.s^2)\end{equation}
    \begin{equation}\rho _1(00.c,r,s^1,s^2) = (00.c,r,s^1,s^2)\end{equation}

    It is straightforward to verify that \(\rho _1\) is well-defined, that \(F(r_1(d),\rho _1(b)) = F(d,b)\) as required. Moreover we see that every finite projection of the output of \(\rho _1\) factors over some finite prefix of the input, and that the resulting maps \(2^N \to  2^M\) carry the uniform distribution to the uniform distribution. This implies that \(\rho _1 c^\omega  = c^\omega \).

\(r_1\) "flattens" the first column of \(d\). Now we can apply \(r_1\) again to the remaining columns to flatten the second column, and similarly apply \(\rho _1\) ignoring the first bit of \(c\) and \(s^1\) to obtain a corresponding permutation. Let \(r_N, \rho _N\) be the map an permutation that flatten the first \(N\) columns. Note that by construction each output bit of \(r_N\) is eventually stable for high \(N\), and the same is true for \(\rho _N\). Let \(r_\infty , \rho _\infty \) be the limiting maps.

    Then when the first \(1\) in the part of \(\rho _\infty (b)\) which selects the \(i\) appears before index \(N\), \(F(r_\infty (d),\rho _\infty (b))\) depends only on the first \(N\) rows, hence is equal to \(F(r_N(d),\rho _\infty (b))\), which is then equal to \(F(d,b)\). Since \(\rho _\infty \) is the limit of permutations, it preserves the uniform measure, and so the composite \(F(d,c^\omega ) = F(r_\infty (d), \rho _\infty (c^\omega )) = F(r_\infty (d), c^\omega )\).

    This proves that \(r_\infty \) preserves the sampling map. But clearly \(r_\infty \) factors over the inclusion of \(BT(A)^\omega  \hookrightarrow  BT(A)^{\omega  \times  \omega }\) as the "flat" sequences of sequences, and equally clearly this inclusion preserves the sampling map. This concludes the proof.
  \qedhere
\end{proof}\end{lemma}

\begin{lemma}\label{efr-HYFI}

  Let \(\mathcal {C}\) be an extensive Markov category with countable Kolmogorov products. Let \(\mathsf {Borel} \to  \mathcal {C}\) be a Markov functor which preserves countable coproducts, countable Kolmogorov products and pullbacks along coproduct inclusions. Then:

\begin{enumerate}\item{}
    For every \(n\), the iterated sampling operator \(BT(A)^{\omega ^n} \to  A\) factors over the map \(BT(A)^{\omega ^n} \to  \bar {\Delta }(A)\).
  
  \item{}
    The composite map \(\coprod _n BT(A)^{\omega ^n} \to  \Delta (A)\) is a split epimorphism, and thus the factorization is unique, giving a well-defined sampling map \(s_A: \Delta (A) \to  A\) for every standard Borel space \(A\).
  
  \item{}
    This sampling map is a \(\omega \)-midpoint homomorphism, in the sense that \(s_A(\sum _i 1/2^i \mu _i) = M_A(s_A \mu _1, s_A \mu _2, \dots )\).
  
  \item{}
    These maps assemble into a functor \(\mathsf {Borel}_{\bar{Delta} } \to  \mathcal {C}\) which extends the unique \(\mathsf {Borel} \to  \mathcal {C}\).
  \end{enumerate}\begin{proof}

    By \cref{efr-X0VA}, we see that we can immediately flatten any nested sequence in \(BT(A)^{\omega ^n}\) into a sequence of dyadic binary trees in \(BT(A)^\omega \), without altering either the sampling map or the represented distribution in \(\bar {\Delta (A)}\). Hence it suffices to prove that the sampling map \(BT(A)^\omega  \to  A\) factors over \(\bar {\Delta (A)}\).

    It is clear that the sum map \(BT(A)^\omega  \to  \bar {\Delta }(A)\) is surjective. We will begin by choosing a particular splitting of this map, assigning to each countably supported measure a "normal form" in \(BT(A)^\omega \)

    Fix an arbitrary total ordering on \(A\). If \(A\) is countable there is no difficulty doing this, if \(A = 2^\omega \), we may use the lexical ordering.

    Given a finitary distribution \(\mu \), we compute its normal form as follows:
  
\begin{enumerate}\item{}
      Let \(N_0\) be the smallest natural so that there exists a dyadic distribution \(\nu _0 = \sum _i \frac {k_i}{2^N} a_i\) so that \(\frac {1}{2}\nu _0 < \mu \) - that is, for every element of \(A\), the probability under \(\frac {1}{2}\nu _0\) is \emph{strictly} less than the probability under \(\mu \).
    
    \item{}
      Let \(\nu _0\) be a dyadic distribution as above, choosing the lowest possible \(a_i\) according to the chosen ordering. Represent \(\nu _0\) as a balanced binary tree of minimal depth, with the leaves increasing left to right.
    
    \item{}
      Repeat these steps with \(\mu  - \frac {1}{2}\nu _0\) to construct \(\nu _1, \nu _2, \dots \).
    \end{enumerate}

    It is clear that this defines a measurable map \(\Delta (A) \to  BT(A)^\omega \), and that this is a section of the map \(\sum _i 1/2^i \nu _i\) in the other direction. 
  
\begin{enumerate}\item{}
      Let \(S: BT(A)^\omega  \times  (2^{\omega } \setminus  0) \times  (2^\omega )^\omega  \to  A\) be the sampling map, which uses the first bitstream to choose the index \(i\) from the sequence, then the \(i\)th one of the remaining streams to sample from the given binary tree.
    
    \item{}
      Note that by construction, if \(\nu _1\) is the first element in the normal form of \(\mu \), and \(\lambda _i\) is another sequence of binary trees so that \(\sum _i 1/2^i \lambda _i\) is equal to \(\mu \), then there exists some finite prefix so that \(1/2 \nu _0 < \sum _{i=0}^N 1/2^i \lambda _i\). Then there exists some rearrangement of this prefix \(\lambda _i'\) so that \(\lambda _1' = \nu _1\) (and \(\lambda '_i = \lambda _i\) when \(i>N\)). Let \(\mathrm {nm}_k: BT(A)^\omega  \to  BT(A)^\omega \) be the map which puts the first \(k\) elements into the normal form like this. Note that for each \(k\) there exists some permutation \(\sigma _k : BT(A)^\omega  \times   (2^{\omega } \setminus  0) \times  (2^\omega )^\omega  \to   (2^{\omega } \setminus  0) \times  (2^\omega )^\omega \) which, for each possible sequence, witnesses this rearrangement (in the same way as \cref{efr-X0VA}).
    
    \item{}
      By construction the \(n\)th coordinate of \(\mathrm {nm}_k(a)\) is equal for all \(k > n\), and so there is a well-defined limiting map \(\mathrm {nm}_\infty  : BT(A)^\omega  \to  BT(A)^\omega \), which (by construction) carries a sequence to its normal form.
    
    \item{}
      For \(\sigma _k\), similarly note that the first \(n\) bits (deciding the index into the sequence) are fixed for \(k > n\), and the same is true for the \(n\) first streams used to sample from the given binary trees.
    
    \item{}
      Let \(\mathrm {nm}_\infty \) and \(\sigma _\infty \) be the limiting maps.
    \end{enumerate}

    Now observe that \(F(a,b) = F(\mathrm {nm}_\infty (a),\sigma _\infty (a,b))\), and that \(\sigma _\infty (a,b)\) for each \(a\) preseves the uniform measure. Hence just as before, this proves that the sampling map commutes with taking the normal form, as we wanted.

    Observe that (for any monad), the Kleisli category is presented by adding morphisms \(s_A: TA \to  A\) for each \(A\), subject to the equations \(s_A\eta _A = 1_A\) and \(s_A s_{TA} = s_A \mu _A\). We have just constructed the \(s_A\), so we now have to prove that they satisfy these equations. Unitality is clear. The multiplicativity condition follows by observing that any finitary distribution on \(\Delta (A)\) can be decomposed as an iterated sum coming from \(\Delta (A)^{\omega ^N}\), at which point this reduces to the mentioned \(\omega \)-midpoint homomorphism property.
  \qedhere
\end{proof}\end{lemma}

    We are now ready for the final step, to construct the extension of the functor \(\mathsf {Borel}_{\bar {\Delta }} \to  \mathcal {C}\) over \(\mathsf {Borel}_{\bar {\Delta }} \hookrightarrow  \mathsf {BorelStoch}\)

\begin{lemma}\label{efr-JK8U}

  Let \(\phi : A \to  (2^\omega )^k\) be a kernel in \(\mathsf {BorelStoch}\), and let \(\vee ^k : (2^\omega )^k \to  2^k\) be the map which is \(1\) in each coordinate if the corresponding sequence contains at least one \(1\) (and zero else). Let \(\mathsf {Borel} \to  \mathcal {C}\) be as above, and let \(F: \mathsf {Borel}_{\bar {\Delta }} \to  \mathcal {C}\) be the unique extensions of \cref{efr-HYFI}. Write \(\bar {F}(\phi )\) for the unique map \(F(A) \to  F((2^\omega )^k)\) given as the limit of \(F\) applied to the finite truncations. Then \(F(\vee ^k)\bar {F}(\phi ) = F(\vee ^k \phi )\)
\begin{proof}

    We proceed by induction on \(k\), the case \(k=0\) being trivial.

    First, consider a kernel \(\alpha : A \to  (2^\omega )^k\) which, on the first \(N\) bits in the first coordinate, is concentrated on a given string \(a = a_1,\dots  a_N\) which contains at least one \(1\). This kernel factors over the inclusion of this component:
    \begin{equation}(2^\omega )^k = \{a\} \otimes  2^{\{N+1, \dots \}} \otimes  (2^\omega )^{k-1} \sqcup  (2^N \setminus  \{a\}) \otimes  2^{\{N+1, \dots \}} \otimes  (2^\omega )^{k-1} \end{equation}

    This factorization is preserved by \(F\) (which preserves coproducts). Note that \(\vee ^k\), when restricted to this summand, factors over the projection
    \((2^\omega )^k \to  (2^\omega )^{k-1}\). So we can write \(F(\vee ^k)\bar {F}(\alpha ) = F((1,\vee ^{k-1}) \pi )\bar {F}(\alpha )\). It is clear that \(F(\pi )\bar {F}(\alpha ) = \bar {F}(\pi  \alpha )\), and so by induction this means \(F(\vee ^k)\bar {F}(\alpha ) = F(\vee ^k \alpha )\).

    By an analogous argument (without passing through the finite truncation), if \(\alpha \) is concentrated on the string \(0 \in  2^\omega \) in the first coordinate, again \(F(\vee ^k)\bar {F}(\alpha ) = F(\vee ^k \alpha )\)

    Now consider the generic kernel \(\phi \) of the lemma. Let \(A = A_0 \sqcup  A_1\), where \(A_0\) consists of those \(a\) for which the probability that the first coordinate is zero is \(\geq  1/2\), and \(A_1\) is the complement, consisting of those for which at least one \(1\) has probability \(>1/2\). Note that by \(\sigma \)-continuity, there exists for each \(a \in  A_1\) some \(N_a\) so that there is a probability \(\geq  1/2\) of an \(1\) in the first \(N_a\) coordinates. Note that the least such \(N_a\) is a measurable function of \(a\), and so \(A_1\) decomposes into a sum \(A_1 = A_1^1 \sqcup  A_1^2 \dots \), where on \(A_1^N\) there is a \(\geq  1/2\) probability of an \(1\) in the first \(N\) coordinates.

    Now let \(\mu _1\) be the kernel defined as follows:
    \begin{enumerate}\item{}
        On \(A_0\), it is deterministically \(0\) in the first coordinate, and given by the conditional distribution given this in the last \(k\) coordinates.
      
      \item{}
        For \(a \in  A_1^N\), choose some distribution \(\nu \) on \(2^N\) so that \(\nu  / 2\) is less than the marginals of \(\phi \), and so that \(\nu \) is concentrated away from \(0\) (this is possible by construction of \(A_1^N\)). Then for each \(b \in  2^N\), let \(\xi _b\) be the measure which is concentrated on that bitstring on the first \(N\) bits, and given by the conditional distribution according to \(\phi \) everywhere else. Let \(\mu _1(a) = \sum _{b\in  2^N} \nu (b) \xi _b\)\end{enumerate}

      Observe that by construction, \(\mu _{1/2} \leq  \phi \). Now we can continue the recursion, defining \(\mu _2, \mu _3, \dots \), and we must necessarily have \(\sum _i 1/2^i \mu _i = \phi \). This sum must be preserved by \(\bar {F}\) (where the meaning of this sum in \(\mathcal {C}\) is given by the canonical sampling morphism \(X^\omega  \to  X\), or equivalently by the dual \(I \to  \omega \)). Therefore we have \begin{equation}F(\vee ^k)\bar {F}(\phi ) = F(\vee ^k)(\sum _i 1/2^i \bar {F}(\mu _i)) = \sum _i 1/2^i F(\vee ^k) \bar {F}(\mu _i) = \sum _i 1/2^i F(\vee ^k \mu _i) = F(\vee ^k \phi ),\end{equation} finishing the proof.
    
\qedhere
\end{proof}\end{lemma}

\begin{lemma}\label{efr-D3QE}

  Let \(\mathsf {Borel} \to  \mathcal {C}\) be as above. Then it admits a unique extension to \(\mathsf {BorelStoch} \to  \mathcal {C}\).

\begin{proof}

  Using \cref{efr-HYFI}, there is a unique extension \(F: \mathsf {Borel}_{\bar {\Delta }} \to  \mathcal {C}\).
  There is a faithful (identity on objects) functor \(\mathsf {Borel}_{\bar {\Delta }} \to  \mathsf {BorelStoch}\) - the only case where this is not full is \(\operatorname {\mathrm {Hom}}(A,2^\omega )\).

  It is apparent that we must extend this functor by defining \(\bar {F}(\phi : A \to  2^\omega ) = \lim _n F(\pi _{2^n}\phi )\), using the universal property of \(F(2^\omega ) = F(2)^\omega \). It is apparent that this preserves independent pairings—the only question is whether this extension is actually functorial. By the limit property it suffices to prove that it is functorial for any composable pair
  \begin{equation}A \xrightarrow {\phi } 2^\omega  \xrightarrow {\psi } K\end{equation}
  with \(K\) finite.

  Let a function \(f: 2^\omega  \to  K\) be \emph{good} if \(\bar {F}(\phi ) ; F(f) = F(\phi  ; f)\) for all kernels \(\phi : A \to  2^\omega \) (for all \(A\)).
  Let an algebra of sets \(\mathbb {A} \subseteq  \mathcal {B}(2^\omega )\) (i.e a collection of subsets stable under finite unions and complements) be called good if every \(\mathbb {A}\)-measurable map \(2^\omega  \to  K\) to a finite set is good. Now we claim:

\begin{enumerate}\item{}
    The class \(\mathbb {A}_0\) of sets of the form \(V \times  2^\omega \) for \(V \subseteq  2^N\), \(N\) finite, is a good algebra.
  
  \item{}
    If \(\mathbb {A}\) is a good algebra, let \(\mathbb {A}^+\) be the smallest algebra containing all countable unions of sets in \(\mathbb {A}\). Then \(\mathbb {A}^+\) is again a good algebra.
  
  \item{}
    Any directed union of good algebras \(\mathbb {A}_0 \subseteq  \mathbb {A}_1 \dots  \subseteq  A_\alpha  \subseteq \) is again a good algebra. 
  \end{enumerate}

  First note that, by Zorn's lemma, this straightforwardly implies the full Borel \(\sigma \)-algebra is a good algebra, which in turn concludes the proof.

  Point 3 holds, since both being an algebra and being good are finitary properties (any map to a finite set which is measurable for the union must be measurable at some finite stage).
  Point 1 is a straightforward consequence of the fact that the projections \(2^\omega  \to  2^N\) are good by construction. So we are left with point 2.

  Let \(\mathbb {A}\) be a good algebra. Let us identify sets by their indicators \(2^\omega  \to  2\). Then we can write any set in \(\mathbb {A}^+\) as
  \begin{equation}g(\vee (f_1^i(-)), \vee  f_2^i(-), \dots  , \vee  f_k^i(-)),\end{equation}
  where \(g: 2^k \to  2\) is some function, \(\vee  : 2^\omega  \to  2\) is the indicator of the sequences with at least one \(1\), and \(f_j^i, 0 \leq  j \leq  k, 0 \leq  i < \infty \) are a bunch of \(\mathbb {A}\)-measurable functions.

  By assumption the pairing of all the \(f\)'s, \(2^\omega  \to  (2^\omega )^k\) is good, and so it suffices to show that the mapping \((\vee )^k : (2^\omega )^k \to  2^k\) is good.
  By induction we may suppose this holds for \((\vee )^{k-1}\), since it is trivial for \(k=0\).

  Now let \(\phi : A \to  (2^\omega )^k\) be some kernel. Write \(A = A_0 + A_1\), where \(A_0\) is the subset where the probability of only zeroes in the first stream is at least \(1/2\), and \(A_1\) is the complement, where the probability of at least one \(1\) in the first stream is \(>1/2\). Observe that by \(\sigma \)-continuity, for each \(a\) there must be some \(N\) so that the probability of at least one \(1\) in the first \(N\) elements of the first stream is at least \(1/2\). Decompose \(A_1\) into \(A_1^1 +  A_1^2 +  \cdots ,\) where for \(a \in  A_1^N\) there is at least a \(1/2\) probability of the first \(1\) being before \(N\). Let \(\psi _1 : A \to  (2^\omega )^k\) be a kernel which on \(A_0\) is given by \((0, \otimes  \mu _{2\dots  k})\) where \(\mu _{2,\dots  k}\) is the kernel giving the conditional distribution of the remaining streams if the first stream is \(0\). On \(A_1^N\), it is given by a linear combination \(\sum _{s \in  2^N \setminus  0} \delta _{s} \otimes  \mu _s\), where \(\mu _s\) is the conditional distribution of the rest of the stream and the remaining streams - we choose this linear combination so that \(1/2\) of it is \(\leq \) the actual probability of each of those prefixes.

  Thus we obtain a decomposition of \(\phi \) as \(\sum _i \frac {1}{2^i} \phi _i\) where for each \(\phi _i\), the distribution on the first coordinate is either concentrated on \(0 \in  2^\omega \), or concentrated away from \(0\) on some finite prefix. Each of these satisfy \(F(\vee ^k)\bar {F}(\phi _i) = F(\vee ^k \phi _i)\). Hence, since \(F(\vee ^k)\) must preserve the infinite linear combination, \(\vee ^k\) must be good as desired.

  Now let \(f: 2^\omega  \to  B\) be a kernel, and let \(\phi : A \to  2^\omega \) again be a kernel. We can factor \(f\) as a measurable map \(\bar {f} : 2^\omega  \times  2^\omega  \to  B\), composed with the uniform \(c^\omega  : I \to  2^\omega \). By the above we have \(F(f) = F(\bar {f})\circ (F(c^\omega ) \otimes  2^\omega )\). By monoidality, it follows that \(\bar {F}(f)\bar {F}(\phi ) = F(f\phi )\) as desired.
\qedhere
\end{proof}\end{lemma}

    Now \cref{efr-TVLB} follows trivially from \cref{efr-AFCK}

  Instead of studying Markov categories, one could study probability monads. The relevant question would then have been to give a universal property for the Giry monad. Our theorem implies this immediately:

\begin{corollary}

    Let \(T\) be a bicommutative, affine monad on \(\mathsf {Borel}\) so that \(T(2)\) contains a coinflip, \(T\) preserves the limits \(\lim _{F \subset  J \mathrm { F finite}} \prod _{i \in  F} X_i\) and pullbacks along coproduct inclusions, and the induced midpoint algebras on \(TX\) are all iterable. Then there is a unique monoidal monad homomorphism \(G \to  T\), where \(G\) is the Giry monad.
  
\begin{proof}

      Subject to these assumptions, \(\mathsf {Borel}_T\) is a Markov category which satisfies the hypotheses of \cref{efr-TVLB}. Hence there is a strong monoidal, identity-on-objects functor \(\mathsf {BorelStoch} = \mathsf {Borel}_G \to  \mathsf {Borel}_T\). Such functors correspond bijectively to morphisms of monads \(G \to  T\) (see e.g. \cite{barr-wells-triples-2005}, Theorem 3.6.3), finishing the proof.
    \qedhere
\end{proof}\end{corollary}

\section{Further work}

    Not all interesting measurable spaces are standard Borel. A simple example is the uncountable product \(\prod _{t \in  \mathbb {R}} \mathbb {R}\), classifying general stochastic processes valued in \(\mathbb {R}\) (although note that many important subsets of this, such as the space of \emph{continuous} stochastic processes, are standard Borel). It would be interesting to look for a generalization of this work to a larger class of measurable spaces, or alternatively to describe in some other terms the initial coinflip Markov category which is \(\kappa \)-extensive, Boolean and admits \(\kappa \)-small Kolmogorov products for larger \(\kappa \).

    Chen's work actually provides the deterministic version of this for all \(\kappa \): the initial \(\kappa \)-complete, (\(\kappa \)-) extensive, Boolean category is given by the opposite of the category of \(\kappa \)-generated, \(\kappa \)-complete Boolean algebras and the \(\kappa \)-continuous homomorphisms between them.

    It is not clear how to define a good probability theory on top of this—having continuum-sized unions seems to rule out many important measures like the Lebesgue measure.

    Fritz and Lorenzin (\cite{fritz-lorenzin-abstractmeas-2025}) have recently described a category \(\mathsf {BaireMeas}\) of \emph{Baire measurable spaces} (those arising as the Baire \(\sigma \)-algebra on a compact Hausdorff space \(X\)), as well as an associated Markov category of kernels \(\mathsf {BaireStoch}\), which admits \emph{all} (small) Kolmogorov products, and which contains \(\mathsf {BorelStoch}\) as a full subcategory.

\begin{conjecture}

\(\mathsf {BaireMeas} \hookrightarrow  \mathsf {BaireStoch}\) is the initial Markov functor from \(\mathsf {BaireMeas}\) into a coinflip Markov category which preserves countable coproducts, pullbacks along coproduct inclusions, and all Kolmogorov products.
    
\end{conjecture}

    This seems amenable to the methods of this paper—what is missing is a concrete description of the general objects of \(\mathsf {BaireMeas}\) as limits of finite ones. However, a universal description of \(\mathsf {BaireMeas}\) itself seems elusive, since it is not Boolean—not every measurable injection between such spaces has a measurable image.

    A general measurable map is difficult to represent computationally. It is therefore interesting to consider variations of this theorem which describe more computationally tractable situations. Note that every compactly supported measure on \(\mathbb {R}^k\) is the image of the uniform measure on \(2^\omega \) under a continuous function, so many measures of interest can be represented in this way (if we work with the extended real line \([-\infty , \infty ]\) instead, every probability measure can be represented).

    Following the discussion above, the right question to ask seems to be, given a countably extensive Cartesian category \(\mathcal {C}\) (for example of topological spaces, or \(\omega \)-dcpos) with products up to a given size, for the free coinflip Markov category \(\overline {\mathcal {C}}\) so that the inclusion preserves countable coproducts, pullbacks along their inclusions and the given Kolmogorov products. It would be interesting to compare this with existing models of continuous probability, such as the category \(\mathsf {TychStoch}\) of Tychonoff spaces and continuous kernels introduced in \cite{markov-supports}, and the category of \(\omega  \mathrm {qbs}\)s and kernels introduced in \cite{vakar-kammar-staton-omega-qbs}.
  
\nocite{*}\bibliographystyle{plain}\bibliography{main.bib}
@article
{lorenzin-zanasi-infinite-tensor-2025, title={Approaching the Continuous from the Discrete: an Infinite Tensor Product Construction}, url={http://arxiv.org/abs/2510.14716}, DOI={10.48550/arXiv.2510.14716}, note={arXiv:2510.14716 [math]}, number={arXiv:2510.14716}, publisher={arXiv}, author={Lorenzin, Antonio and Zanasi, Fabio}, year={2025}, month={Oct} }

@article
{fritz-lorenzin-abstractmeas-2025, title={Categories of abstract and noncommutative measurable spaces}, url={http://arxiv.org/abs/2504.13708}, DOI={10.48550/arXiv.2504.13708}, note={arXiv:2504.13708 [math]}, number={arXiv:2504.13708}, publisher={arXiv}, author={Fritz, Tobias and Lorenzin, Antonio}, year={2025}, month={Apr} }

@article
{law-large-nums-fritz-etal, title={Empirical Measures and Strong Laws of Large Numbers in Categorical Probability}, url={http://arxiv.org/abs/2503.21576}, DOI={10.48550/arXiv.2503.21576}, number={arXiv:2503.21576}, publisher={arXiv}, author={Fritz, Tobias and Gonda, Tomas and Lorenzin, Antonio and Perrone, Paolo and Mohammed, Areeb Shah}, year={2025}, month={Mar} }

@article
{markov-supports, title={Absolute continuity, supports and idempotent splitting in categorical probability}, url={http://arxiv.org/abs/2308.00651}, DOI={10.48550/arXiv.2308.00651}, abstractNote={Markov categories have recently turned out to be a powerful high-level framework for probability and statistics. They accommodate purely categorical definitions of notions like conditional probability and almost sure equality, as well as proofs of fundamental results such as the Hewitt-Savage 0/1 Law, the de Finetti Theorem and the Ergodic Decomposition Theorem. In this work, we develop additional relevant notions from probability theory in the setting of Markov categories. This comprises improved versions of previously introduced definitions of absolute continuity and supports, as well as a detailed study of idempotents and idempotent splitting in Markov categories. Our main result on idempotent splitting is that every idempotent measurable Markov kernel between standard Borel spaces splits through another standard Borel space, and we derive this as an instance of a general categorical criterion for idempotent splitting in Markov categories.}, number={arXiv:2308.00651}, publisher={arXiv}, author={Fritz, Tobias and Gonda, Tomas and Lorenzin, Antonio and Perrone, Paolo and Stein, Dario}, year={2023}, month={Sep} }

@article
{fritz-gonda-perrone-rischel-rep, title={Representable Markov categories and comparison of statistical experiments in categorical probability}, volume={961}, ISSN={0304-3975}, DOI={10.1016/j.tcs.2023.113896}, journal={Theoretical Computer Science}, author={Fritz, Tobias and Gonda, Tomas and Perrone, Paolo and Fjeldgren Rischel, Eigil}, year={2023}, month={Jun}, pages={113896} }

@article
{fritz-liang-freemarkov-2023, title={Free gs-monoidal categories and free Markov categories}, volume={31}, ISSN={0927-2852, 1572-9095}, DOI={10.1007/s10485-023-09717-0}, note={arXiv:2204.02284 [math]}, number={2}, journal={Applied Categorical Structures}, author={Fritz, Tobias and Liang, Wendong}, year={2023}, month={Apr}, pages={21} }

@article
{chen-universal-stdborel-2019, title={A Universal Characterization of Standard Borel Spaces}, volume={88}, DOI={10.1017/jsl.2022.85}, number={2}, journal={The Journal of Symbolic Logic}, author={Chen, Ruiyuan}, year={2022}, month={Dec}, pages={510–539} }

@article
{fritz-gonda-perrone-2021, title={De Finetti’s Theorem in Categorical Probability}, volume={2}, ISSN={2689-6931}, url={http://arxiv.org/abs/2105.02639}, DOI={10.31390/josa.2.4.06}, note={arXiv:2105.02639 [math]}, number={4}, journal={Journal of Stochastic Analysis}, author={Fritz, Tobias and Gonda, Tomáš and Perrone, Paolo}, year={2021}, month={Nov} }

@article
{rischel-fritz-infinite-products, title={Infinite products and zero-one laws in categorical probability}, volume={2}, ISSN={2631-4444}, DOI={10.32408/compositionality-2-3}, journal={Compositionality}, author={Fritz, Tobias and Rischel, Eigil Fjeldgren}, year={2020}, month={Aug}, pages={3},language={en} }

@article
{fritz-synthetic-markov-cats, title={A synthetic approach to Markov kernels, conditional independence and theorems on sufficient statistics}, volume={370}, ISSN={00018708}, DOI={10.1016/j.aim.2020.107239}, abstractNote={We develop Markov categories as a framework for synthetic probability and statistics, following work of Golubtsov as well as Cho and Jacobs. This means that we treat the following concepts in purely abstract categorical terms: conditioning and disintegration; various versions of conditional independence and its standard properties; conditional products; almost surely; suï¬ƒcient statistics; as well as versions of theorems on suï¬ƒcient statistics due to Fisher-Neyman, Basu, and Bahadur. Besides the conceptual clarity oï¬€ered by our categorical setup, its main advantage is that it provides a uniform treatment of various types of probability theory, including discrete probability theory, measure-theoretic probability with general measurable spaces, Gaussian probability, Markov processes of either of these kinds, and many others.}, url={https://arxiv.org/abs/1908.07021}, journal={Advances in Mathematics}, author={Fritz, Tobias}, year={2020}, month={Aug}, pages={107239},language={en} }

@article
{vakar-kammar-staton-omega-qbs, title={A Domain Theory for Statistical Probabilistic Programming}, volume={3}, ISSN={2475-1421}, DOI={10.1145/3290349}, note={arXiv:1811.04196 [cs]}, number={POPL}, journal={Proceedings of the ACM on Programming Languages}, author={Vákár, Matthijs and Kammar, Ohad and Staton, Sam}, year={2019}, month={Jan}, pages={1–29} }

@article
{fritz-stochmat-2009, title={A presentation of the category of stochastic matrices}, rights={arXiv.org perpetual, non-exclusive license}, url={https://arxiv.org/abs/0902.2554}, DOI={10.48550/ARXIV.0902.2554}, publisher={arXiv}, author={Fritz, Tobias}, year={2009} }

@article
{golubtsov-kleisli, title={Monoidal Kleisli Category as a Background for Information Transformers Theory}, abstractNote={We consider any uniform class of information transformers (ITs) as a family of morphisms of a monoidal category that contains a subcategory (of deterministic ITs) with finite products and satisfies certain set of axioms. Besides, many IT-categories can be constructed as Kleisli categories. The ingredients for this construction are: a base category (of deterministic ITs); a functor, producing objects of â€œdistributionsâ€; a natural transformation, representing â€œindependent product of distributionsâ€. The paper also generalizes Bayesian approach to decision-making problems and studies informativeness of ITs. It shows that classes of equivalent ITs form a partially ordered bounded Abelian monoid. Several examples of concrete IT-categories are examined.}, author={Golubtsov, Peter V}, year={2002}, pages={24},language={en} }

@book
{johnstone-elephant-vol1, address={Oxford, England}, title={Sketches of an Elephant: A Topos Theory Compendium, Volume 1}, publisher={Clarendon Press}, author={Johnstone, Peter T.}, year={2002} }

@inproceedings
{escardo-simpson-universal-interval-2001, title={A universal characterization of the closed Euclidean interval}, ISSN={1043-6871}, url={https://ieeexplore.ieee.org/document/932488/}, DOI={10.1109/LICS.2001.932488}, booktitle={Proceedings 16th Annual IEEE Symposium on Logic in Computer Science}, author={Escardo, M.H. and Simpson, A.K.}, year={2001}, month={June}, pages={115–125} }

@article
{carboni-lack-walters-extensive, title={Introduction to extensive and distributive categories}, volume={84}, ISSN={0022-4049}, DOI={10.1016/0022-4049(93)90035-R}, number={2}, journal={Journal of Pure and Applied Algebra}, author={Carboni, Aurelio and Lack, Stephen and Walters, R. F. C.}, year={1993}, month={Feb}, pages={145–158} }

@inproceedings
{giry-1982, address={Berlin, Heidelberg}, title={A categorical approach to probability theory}, ISBN={978-3-540-39041-1}, DOI={10.1007/BFb0092872}, booktitle={Categorical Aspects of Topology and Analysis}, publisher={Springer}, author={Giry, Michèle}, editor={Banaschewski, B.}, year={1982}, pages={68–85},language={en} }

@article
{barr-wells-triples-2005, title={TOPOSES, TRIPLES AND THEORIES}, author={Barr, Michael and Wells, Charles}, pages={302},language={en} }
\end{document}